\documentclass[twoside,12pt]{article}
\usepackage{amssymb,amsmath,bm,euscript}
\usepackage{graphicx,color,caption}

\usepackage[colorlinks=true]{hyperref}

\usepackage{algorithm,algpseudocode}
\usepackage{braket}

\newtheorem{proposition}{Proposition}[section]
\newtheorem{theorem}[proposition]{Theorem}

\newtheorem{corollary}[proposition]{Corollary}
\newtheorem{definition}[proposition]{Definition}

\newcommand{\qed}{\hphantom{.}\hfill $\Box$\medbreak}

\def\T{\mathcal{T}}

\def\A{{\mathcal{A}}}

\def\B{{\mathcal{B}}}

\def\T{\mathcal{T}}

\def\x{{\bf x}}
\def\y{{\bf y}}
\def\z{{\bf z}}
\def\w{{\bf w}}
\def\s{{\bf s}}
\def\uu{{\bf u}}
\def\vv{{\bf v}}
\def\0{{\bf 0}}

\title{\bf{Spectral Norm and Nuclear Norm of a Third Order Tensor}}
\author{ \hspace{1mm} Liqun Qi\thanks{Department of Applied
    Mathematics, The Hong Kong Polytechnic University, Hung Hom,
    Kowloon, Hong Kong; ({\tt liqun.qi@polyu.edu.hk}). This author's work was supported by the Hong Kong
    Research Grant Council (Grant No.  PolyU 15300715, 15301716 and 15300717). }
    \ and\
    Shenglong Hu\thanks{Department of Mathematics, School of Science, Hangzhou Dianzi University, Hangzhou 310018 China; ({\tt shenglonghu@hdu.edu.cn}). This author's work was supported by NSFC (Grant No.  11771328).}}

\begin{document}
\date{\today}
\maketitle

\begin{abstract}
The spectral norm and the nuclear norm of a third order tensor play an important role in the tensor completion and recovery problem.   We show that the spectral norm of a third order tensor is equal to the square root of the spectral norm of three positive semi-definite biquadratic tensors, and the square roots of the nuclear norms of those three  positive semi-definite biquadratic tensors are lower bounds of the nuclear norm of that third order tensor.  This provides a way to estimate and to evaluate the spectral norm and the nuclear norm of that third order tensor.  Some upper and lower bounds for the spectral norm and nuclear norm of a third order tensor, by spectral radii and nuclear norms of some symmetric matrices, are presented.


\end{abstract}


\section{Introduction}

The spectral norm and the nuclear norm of a third order tensor play an important role in the tensor completion and recovery problem \cite{SGCH19, YZ16}.  It is NP-hard to compute them \cite{FL17}.  It is an active research topic to study them more \cite{Hu15, JYZ17, Li16}.

In this paper, unless otherwise stated, all the discussions will be carried out in the filed of real numbers.
The spectral norm of a third norm is the largest singular value of that tensor.   The nuclear norm is the dual norm of the spectral norm.   Hence,
singular values of a third order tensor form the base of the spectral norm and the nuclear norm.  Recall that the product of a (maybe rectangular) matrix and the transpose of that matrix is a positive semi-definite symmetric (square) matrix.   There is a one to one equality  between the singular values of the original matrix and the square roots of the eigenvalues of that positive semi-definite symmetric matrix.  Then the spectral norm of the original  matrix is equal to the square root of the spectral radius of that positive semi-definite symmetric matrix.   Does such a relation still exist for a third order tensor?   In the next section, we give a firm answer to this question.  We show that if we make contraction of a third order tensor with itself on one index, then we get a  positive semi-definite biquadratic tensor.   A real number is a singular value of that third order tensor if and only if it is the square root of an M-eigenvalue of that positive semi-definite biquadratic tensor.   Thus, the spectral norm of that third order tensor is the square root of the spectral norm of that  positive semi-definite biquadratic tensor.

In Section 3, we show that the square root of the nuclear norm of that  positive semi-definite biquadratic tensor is a lower bound of the nuclear norm of that third order tensor.   The equality may not hold in general.

The equality between the spectral norm of a third order tensor and the spectral norm of a  positive semi-definite biquadratic tensor does not change the complexity of the problem, but provides us an alternative way to attack the problem.  In Sections 4 and 5, by this relation, we present several upper and lower bounds for the spectral norm of a third order tensor, by spectral radii of some symmetric matrices.   In Section 6, we establish some relations between these upper and lower bounds, and thus give a range for the spectral norm of that third order tensors.

In Section 7, we present some lower bounds for the nuclear norm of a third order tensor, by the nuclear norms of some symmetric matrices.

Some final remarks are made in Section 8.

\section{Spectral Norm}

Suppose that $d_1, d_2$ and $d_3$ are positive integers.  Without loss of generality, we may assume that $d_1 \le d_2 \le d_3$.

Let $\Re^{d_1 \times d_2 \times d_3}$ be the space of third order tensors of dimension $d_1 \times d_2 \times d_3$.   The singular values of a tensor $\A = (a_{ijk}) \in \Re^{d_1 \times d_2 \times d_3}$ are defined as follows \cite{Lim05}.

\begin{definition}
A real number $\lambda$ is called a singular value of $\A$ if there are vectors $\x = (x_1, \cdots, x_{d_1})^\top \in \Re^{d_1}, \y = (y_1, \cdots, y_{d_2})^\top \in \Re^d_2, \z = (z_1, \cdots, z_{d_3})^\top \in \Re^{d_3}$ such that the following equations are satisfied:
For $i = 1, \cdots, d_1$,
\begin{equation} \label{e1}
\sum_{j=1}^{d_2}\sum_{k=1}^{d_3} a_{ijk}y_jz_k = \lambda x_i;
\end{equation}
For $j = 1, \cdots, d_2$,
\begin{equation} \label{e2}
\sum_{i=1}^{d_1}\sum_{k=1}^{d_3} a_{ijk}x_iz_k = \lambda y_j;
\end{equation}
For $k = 1, \cdots, d_3$,
\begin{equation} \label{e3}
\sum_{i=1}^{d_1}\sum_{j=1}^{d_2} a_{ijk}x_iy_j = \lambda z_k;
\end{equation}
and
\begin{equation} \label{e4}
\x^\top \x = \y^\top \y = \z^\top \z = 1.
\end{equation}
Then $\x, \y$ and $\z$ are called the corresponding singular vectors.
\end{definition}

If $\lambda$ is a singular value of $\A$, with singular vectors $\x, \y$ and $\z$, then by definition, $-\lambda$ is also a singular value of $\A$, with singular vector $-\x, -\y$ and $-\z$.
For $\A =  (a_{ijk}), \B = (b_{ijk}) \in \Re^{d_1 \times d_2 \times d_3}$, their inner product is defined as
$$
\langle \A, \B \rangle := \sum_{i=1}^{d_1} \sum_{j=1}^{d_2} \sum_{k=1}^{d_3} a_{ijk}b_{ijk}.
$$
In a special case, if $\B$ is rank-one, i.e., $\B = (b_{ijk}) = \x \otimes \y \otimes \z$ for some nonzero vectors $\x \in \Re^{d_1}, \y \in \Re^{d_2}, \z \in \Re^{d_3}$, or equivalently $b_{ijk} = x_iy_jz_k$ for $i = 1, \cdots, d_1, j = 1, \cdots, d_2$ and $k = 1, \cdots, d_3$, then
$$\langle \A, \x \otimes \y \otimes \z \rangle \equiv \sum_{i=1}^{d_1} \sum_{j=1}^{d_2} \sum_{k=1}^{d_3} a_{ijk}x_iy_jz_k.$$

\begin{definition}
The spectral norm of $\A \in \Re^{d_1 \times d_2 \times d_3}$ is defined \cite{FL17, Hu15, JYZ17, Li16} as
\begin{equation} \label{n1}
\| \A \| : = \max \left\{ \langle \A, \x \otimes \y \otimes \z \rangle  : \x^\top \x = \y^\top \y = \z^\top \z = 1, \x \in \Re^{d_1}, \y \in \Re^{d_2}, \z \in \Re^{d_3} \right\}.
\end{equation}
\end{definition}
Then the  spectral norm of $\A$ is equal to the largest singular value of $\A$ \cite{FL17, Hu15, JYZ17, Li16}.

\medskip
We now consider  biquadratic tensors.

\begin{definition}
Let $\Re^{d_1 \times d_2 \times d_1 \times d_2}$ be the space of fourth order tensors of dimension $d_1 \times d_2 \times d_1 \times d_2$.   Let $\T = (t_{ijpq}) \in \Re^{d_1 \times d_2 \times d_1 \times d_2}$.  The tensor $\T$ is called biquadratic if for all $i, p = 1, \cdots, d_1$ and $j, q = 1, \cdots, d_2$, we have
$$t_{ijpq} = t_{pjiq} = t_{pqij}.$$
The tensor $\T$ is called positive semi-definite if for any $\x \in \Re^{d_1}$ and $\y \in \Re^{d_2}$,
$$\langle \T, \x \otimes \y \otimes \x \otimes \y \rangle \equiv \sum_{i, p =1}^{d_1} \sum_{j, q = 1}^{d_2} t_{ijpq}x_iy_jx_py_q \ge 0.$$
The tensor $\T$ is called positive definite if for any $\x \in \Re^{d_1}, \x^\top \x = 1$ and $\y \in \Re^{d_2}, \y^\top \y = 1$,
$$\langle \T, \x \otimes \y \otimes \x \otimes \y \rangle \equiv \sum_{i, p =1}^{d_1} \sum_{j, q = 1}^{d_2} t_{ijpq}x_iy_jx_py_q > 0.$$
The spectral norm of $\T$ is defined by
\begin{equation} \label{n3}
\| \T \| := \max \left\{ \left| \langle \T, \x \otimes \y \otimes \x \otimes \y \rangle \right| :  \x^\top \x = \y^\top \y = 1, \x \in \Re^{d_1}, \y \in \Re^{d_2} \right\}.
\end{equation}
\end{definition}
We may check that $\| \cdot \|$ defines a norm in $\Re^{d_1 \times d_2 \times d_1 \times d_2}$.

\begin{definition}
Suppose that $\T = (t_{ijpq}) \in \Re^{d_1 \times d_2 \times d_1 \times d_2}$  is biquadratic.   A number $\mu$ is called an M-eigenvalue of $\T$ if there are vectors  $\x = (x_1, \cdots, x_{d_1})^\top \in \Re^{d_1}, \y = (y_1, \cdots, y_{d_2})^\top \in \Re^d_2$ such that the following equations are satisfied:
For $i = 1, \cdots, d_1$,
\begin{equation} \label{e5}
\sum_{p=1}^{d_1}\sum_{j, q=1}^{d_2} t_{ijpq}y_jx_py_q = \mu x_i;
\end{equation}
For $j = 1, \cdots, d_2$,
\begin{equation} \label{e6}
\sum_{i,p=1}^{d_1}\sum_{q=1}^{d_2} t_{ijpq}x_ix_py_q = \mu y_j;
\end{equation}
and
\begin{equation} \label{e7}
\x^\top \x = \y^\top \y = 1.
\end{equation}
Then $\x$ and $y$ are called the corresponding M-eigenvectors.
\end{definition}

\begin{theorem} \label{t1}
Suppose that $\T = (t_{ijpq}) \in \Re^{d_1 \times d_2 \times d_1 \times d_2}$  is biquadratic.  Then its M-eigenvalues always exist.  The spectral norm of $\T$ is equal to the largest absolute value of its M-eigenvalues.   Furthermore, $\T$ is positive semi-definite if and only if all of its M-eigenvalues are nonnegative; $\T$ is positive definite if and only if all of its M-eigenvalues are positive.  If $\T$ is positive semi-definite, then its spectral norm is equal to its largest M-eigenvalue.
\end{theorem}
{\bf Proof} Consider the optimization problem
\begin{equation} \label{e8}
\min \left\{ \langle \T, \x \otimes \y \otimes \x \otimes \y \rangle :  \x^\top \x = \y^\top \y = 1, \x \in \Re^{d_1}, \y \in \Re^{d_2} \right\}.
\end{equation}
Since the objective function is continuous and the feasible region is compact, this optimization problem always has an optimal solution.  Since the linear independence constraint qualification in optimization is satisfied, the optimality condition holds at that optimal solution.   By optimization theory, the optimality condition of (\ref{e8}) has the form (\ref{e5}-\ref{e7}), and the optimal Langrangian multiplier $\mu$ always exists at the solution.  This shows that $\T$ always has an M-eigenvalue.

Suppose that $\mu$ is an M-eigenvalue of $\T$ with corresponding vectors $\x$ and $\y$.  By (\ref{e5}) and (\ref{e6}), we have
$$\mu = \langle \T, \x \otimes \y \otimes \x \otimes \y \rangle.$$
By this and (\ref{n3}), the spectral norm of $\T$ is equal to the largest absolute value of its M-eigenvalues.
By this and (\ref{e8}), $\T$ is positive semi-definite if and only if all of its M-eigenvalues are nonnegative; $\T$ is positive definite if and only if all of its M-eigenvalues are positive. If $\T$ is positive semi-definite, then all of its M-eigenvalues are nonnegative.  This implies that its spectral norm is equal to its largest M-eigenvalue in this case.
\qed

For $d_1 = d_2 =3$, the elastic tensor in solid mechanics falls in the form of $\T$, with two additional symmetric properties between indices $i$ and $j$, and between indices $p$ and $q$.   Then, the positive definiteness condition of $\T$ corresponds the strong ellipticity condition in solid mechanics.  In 2009, M-eigenvalues were introduced for the elastic tensor to characterize the strong ellipticity condition in \cite{QDH09}.   An algorithm for computing the largest M-eigenvalue was presented in \cite{WQZ09}.  Also see \cite{QCC18} for details.  Here, we extend M-eigenvalues to general biquadratic tensors and study their spectral norms.

\medskip

For $\A =  (a_{ijk}) \in \Re^{d_1 \times d_2 \times d_3}$, consider its contraction with itself on the third index, $\T^{(3)} = \left(t^{(3)}_{ijpq}\right) \in \Re^{d_1 \times d_2 \times d_1 \times d_2}$, defined by
\begin{equation} \label{e9}
t^{(3)}_{ijpq} = \sum_{k=1}^{d_3} a_{ijk}a_{pqk}.
\end{equation}
Then $\T^{(3)}$ is biquadratic.   For any $\x \in \Re^{d_1}$ and $\y \in \Re^{d_2}$,
$$\langle \T^{(3)}, \x \otimes \y \otimes \x \otimes \y \rangle = \sum_{k=1}^{d_3} \left( \sum_{i=1}^{d_1} \sum_{j= 1}^{d_2} a_{ijk}x_iy_j \right)^2 \ge 0.$$
Hence $\T^{(3)}$ is also positive semi-definite.

\begin{theorem} \label{t2}
Let $\A = (a_{ijk}) \in \Re^{d_1 \times d_2 \times d_3}$ and $\T^{(3)} = \left(t^{(3)}_{ijpq}\right) \in \Re^{d_1 \times d_2 \times d_1 \times d_2}$ be constructed as above.   Then $\lambda$ is a nonzero singular value of $\A$, with $\x \in \Re^{d_1}, \y \in \Re^{d_2}$ and $\z \in \Re^{d_3}$ as its corresponding singular vectors, if and only if it is a square root of an M-eigenvalue of $\T^{(3)}$, with $\x$ and $\y$ as corresponding M-eigenvectors.   This also implies that the spectral norm of $\A$ is equal to the square root of the largest M-eigenvalue of $\T^{(3)}$.
\end{theorem}
{\bf Proof} Suppose that $\lambda \not = 0$ is a singular value of $\A$, with corresponding singular vectors $\x, \y$ and $\z$, satisfying (\ref{e1}-\ref{e4}).  Multiplying (\ref{e1}) and (\ref{e2}) by $\lambda$ and substituting
$$\lambda z_k = \sum_{p=1}^{d_1} \sum_{q=1}^{d_2} a_{pqk}x_py_q$$
into these two equations, we see that $\mu = \lambda^2$ is an M-eigenvalue of $\T^{(3)}$, with $\x$ and $\y$ as the corresponding M-eigenvectors.

On the other hand, assume that $\mu = \lambda^2 \not = 0$ is an M-eigenvalue of $\T^{(3)}$, with corresponding M-eigenvectors $\x$ and $\y$, satisfying (\ref{e5}-\ref{e7}), where $\T^{(3)}$ is constructed as above.   Let
$\z = (z_1, \cdots, z_{d_3})^\top$ with
$$ z_k = {1 \over \lambda} \sum_{i=1}^{d_1} \sum_{j=1}^{d_2}a_{ijk}x_iy_j.$$
Then (\ref{e3}) is satisfied.
$$\begin{aligned}
\z^\top \z & =  {1 \over \lambda^2} \sum_{k=1}^{d_3} \left(\sum_{i=1}^{d_1} \sum_{j=1}^{d_2}a_{ijk}x_iy_j \sum_{p=1}^{d_1} \sum_{q=1}^{d_2}a_{pqk}x_py_q \right) \\
& = {1 \over \mu}\sum_{i, p=1}^{d_1}\sum_{j, q=1}^{d_2} \left( \sum_{k=1}^{d_3} a_{ijk}a_{pqk}\right)x_iy_jx_py_q \\
& = {1 \over \mu} \sum_{i=1}^{d_1} \left( \sum_{p=1}^{d_1}\sum_{j, q=1}^{d_2} t^{(3)}_{ijpq}y_jx_py_q \right)x_i \\
& = \sum_{i=1}^{d_1} x_i^2 \\
& = 1.
\end{aligned}$$
This proves (\ref{e4}).
We also have
$$\begin{aligned}
\sum_{j=1}^{d_2}\sum_{k=1}^{d_3} a_{ijk}y_jz_k & = {1 \over \lambda} \sum_{p=1}^{d_1} \sum_{j, q=1}^{d_2}\sum_{k=1}^{d_3} a_{ijk}a_{pqk}x_py_jy_q\\
& = {1 \over \lambda} \sum_{p=1}^{d_1} \sum_{j, q=1}^{d_2} t^{(3)}_{ijpq}x_py_jy_q\\
& = {\mu x_i \over \lambda}\\
& = \lambda x_i.
\end{aligned}
$$
This proves (\ref{e1}).  We may prove (\ref{e2}) similarly.  Hence, $\lambda$ is a singular value of $\A$, with $\x, \y$ and $\z$ as the corresponding singular vectors.

By Theorem \ref{t1}, we now conclude that the spectral norm of $\A$ is equal to the square root of the largest M-eigenvalue of $\T^{(3)}$.
\qed

{\bf Example 1} Let the entries of $\A = (a_{ijk}) \in \Re^{2\times 2 \times 3}$ be
$$ \begin{aligned} a_{111} & =& 4, \ a_{121} & =& 1,  \ a_{112} & =& 3, \ a_{122} & = & 2, \ a_{113} & = & 2, \ a_{123} & = & -1, \\ a_{211} & =& -1, \ a_{221} & =& 2, \ a_{212} & =& -5,  \ a_{222} & =& 1, \ a_{213} & =& 3, \ a_{223} & =& 4. \end{aligned}$$
Calculate the spectral norm of $\A$ by definition, we see that the spectral norm of $\A$ is $6.7673$.
Then the entries of $\T^{(3)} =\left(t^{(3)}_{ijpq}\right)$ are $t^{(3)}_{1111} = 29$, $t^{(3)}_{1112} = t^{(3)}_{1211} = 8$, $t^{(3)}_{1121} = t^{(3)}_{2111} = -13$, $t^{(3)}_{1212} = 6$, $t^{(3)}_{1221} = t^{(3)}_{2112} = -14$, $t^{(3)}_{1122} = t^{(3)}_{2211} = 19$, $t^{(3)}_{2121} = 35$, $t^{(3)}_{1222} = t^{(3)}_{2212} = 0$, $t^{(3)}_{2122} = t^{(3)}_{2221} = 5$, $t^{(3)}_{2222} = 21$.
Calculate the spectral norm of $\T^{(3)}$ by definition, we see that the spectral norm of $\T$ is $45.7959$.   Its square root is $6.7673$, which is equal to the spectral norm of $\A$.
\qed

\begin{corollary}  \label{c1}
We may also consider the contraction of $\A$ and itself over its second index or the first index.   Then we have a tensor $\T^{(2)}$ in $\Re^{d_1 \times d_3 \times d_1 \times d_3}$ and a tensor $\T^{(1)}$ in $\Re^{d_2 \times d_3 \times d_2 \times d_3}$.  Theorem \ref{t2} is true for $\A$ and these two  positive semi-definite biquadratic tensors $\T^{(2)}$ and $\T^{(1)}$ too.
\end{corollary}

Our numerical computation confirms the results of Theorem \ref{t2} and Corollary \ref{c1}.

\section{Nuclear Norm}

The nuclear norm is somewhat more important in the tensor completion and recovery problem \cite{SGCH19, YZ16}.

\begin{definition}
The nuclear norm of $\A \in \Re^{d_1 \times d_2 \times d_3}$ is defined \cite{FL17, Li16} as
\begin{equation} \label{n2}
\|\A \|_* := \inf \left\{ \sum_{i=1}^r |\lambda_i| :  \A = \sum_{i=1}^r \lambda_i \uu_i \otimes \vv_i \otimes \w_i, {\uu_i^\top \uu_i = \vv_i^\top \vv_i = \w_i^\top \w_i = 1, \atop \lambda_i\in \Re, \uu_i \in \Re^{d_1}, \vv_i \in \Re^{d_2}, \w_i \in \Re^{d_3},} i=1, \cdots, r \right\}.
\end{equation}
\end{definition}
Then we have \cite{FL17, Li16}
\begin{equation}\label{eq:dual}
\|\A \|_* := \max \left\{ \langle \A, \B \rangle : \| \B \| = 1, \B \in \Re^{d_1 \times d_2 \times d_3} \right\}.
\end{equation}

We may define the nuclear norm of a tensor in $\Re^{d_1 \times d_2 \times d_1 \times d_2}$ similarly.

\begin{definition}
The nuclear norm of $\T \in \Re^{d_1 \times d_2 \times d_1 \times d_2}$ is defined  as
\begin{equation} \label{n4}
\|\T \|_* := \inf \left\{ \sum_{i=1}^r |\lambda_i| :  \T = \sum_{i=1}^r \lambda_i \uu_i \otimes \vv_i \otimes \w_i \otimes \s_i, {\uu_i^\top \uu_i = \vv_i^\top \vv_i = \w_i^\top \w_i = \s_i^\top \s_i = 1, \atop \lambda_i\in\Re, \uu_i, \w_i \in \Re^{d_1}, \vv_i, \s_i \in \Re^{d_2},} i=1, \cdots, r \right\}.
\end{equation}
\end{definition}
 Then we have the following theorem.
 \begin{theorem} \label{t3}
 Suppose that  $\A = (a_{ijk}) \in \Re^{d_1 \times d_2 \times d_3}$, and $\T^{(3)} = \left(t^{(3)}_{ijpq}\right)$ is constructed by (\ref{e9}). Assume $\| \A \|_*$ and $\|\T^{(3)}\|_*$ are defined by (\ref{n2}) and (\ref{n4}) respectively.   Then
 \begin {equation} \label{e14}
 \| \A \|_*^2 \ge \left\| \T^{(3)} \right\|_*\ge \frac{1}{d_3}\| \A\|_*^2.
 \end{equation}
 \end{theorem}
 {\bf Proof} For any $\epsilon > 0$, by (\ref{n2}), we have positive integer $r$ and $\uu_i \in \Re^{d_1}, \vv_i \in \Re^{d_2}, \w_i \in \Re^{d_3}$ such that
 $$\uu_i^\top \uu_i = \vv_i^\top \vv_i = \w_i^\top \w_i = 1,$$
 for $i = 1, \cdots, r$, and
 $$\A = \sum_{i=1}^r \lambda_i \uu_i \otimes \vv_i \otimes \w_i$$
 and
 $$\|\A \|_* + \epsilon \ge \sum_{i=1}^r |\lambda_i|.$$
 By (\ref{e9}), we have
 $$\T^{(3)} = \sum_{i, j=1}^r \lambda_i\lambda_j\alpha_{ij} \uu_i \otimes \vv_i \otimes \uu_j \otimes \vv_j,$$
 where $\alpha_{ij}=\w_i^\top \w_j$.
 Then by (\ref{n4}), we have
 $$\left(\|\A \|_* + \epsilon\right)^2 \ge \|\T^{(3)}\|_*$$
 for any $\epsilon > 0$.
 This proves the first inequality in (\ref{e14}).

 For the lower bound in \eqref{e14}, suppose that $\B\in\Re^{d_1\times d_2\times d_3}$ is such that
 \[
 \|\B\|=1\ \text{and }\langle\A,\B\rangle =\|\A\|_*.
 \]
 For simplicity of notation, denote by the $d_1\times d_2$ matrix $[a_{\cdot\cdot k}]$ as $A_k$ for all $k=1,\dots,d_3$. Similarly, we have $d_3$ matrices $B_k$'s for $\B$. Since $\|\A\|_*$ is the maximum of $\langle\A,\B\rangle$ over all tensors $\B$ with unit spectral norm, and the spectral norm is defined by maximizing a multilinear function over the joint sphere (cf.\ \eqref{n1}), we must have that
 \[
 \langle A_k,B_k\rangle\geq 0\ \text{for all }k=1,\dots,d_3\ \text{and }\|\A\|_*=\sum_{k=1}^{d_3} \langle A_k,B_k\rangle.
 \]
 Let the tensor $\mathcal S$ be defined similarly to $\T^{(3)}$ for $\A$, i.e., $\mathcal S=\sum_{k=1}^{d_3}B_k\otimes B_k$. It follows from Theorems \ref{t1} and \ref{t2} that
 \[
 \|\mathcal S\|=1.
 \]
 Then, by \eqref{eq:dual}, we have
 \[
 \left\|\T^{(3)}\right\|_*\geq \langle\T^{(3)},\mathcal S\rangle=\sum_{k=1}^{d_3}\langle A_k,B_k\rangle^2\geq \frac{1}{d_3}\left(\sum_{k=1}^{d_3}\langle A_k,B_k\rangle\right)^2= \frac{1}{d_3}\|\A\|_*^2.
 \]
 The second inequality in \eqref{e14} is thus proved.
  \qed

Numerical computations show that strict inequality may hold in (\ref{e14}).

\begin{corollary} \label{c2}
We may also consider the contraction of $\A$ and itself over its second index or the first index.   Then we have a tensor $\T^{(2)}$ in $\Re^{d_1 \times d_3 \times d_1 \times d_3}$ and a tensor $\T^{(1)}$ in $\Re^{d_2 \times d_3 \times d_2 \times d_3}$.  Theorem \ref{t3} is true for $\A$ and these two  positive semi-definite biquadratic tensors $\T^{(2)}$ and $\T^{(1)}$ too.
\end{corollary}

Numerical computation shows that the nuclear norms of these three  positive semi-definite biquadratic tensors can be different for a third order tensor $\A$.

\section{Upper Bounds}

Theorems \ref{t2} and \ref{t3} connect the spectral norm and nuclear norm of a third order tensor with the spectral norms and nuclear norms of three  positive semi-definite biquadratic tensors.  This does not change the complexity of the problem.  But they provide us an alternative way to attack the problem.   In particular, a  biquadratic tensor has more structure such as the diagonal structure.  In 2009, Wang, Qi and Zhang \cite{WQZ09} presented a practical method for the largest M-eigenvalue of a  biquadratic tensor.
Thus, we may apply that method to compute the spectral norm of a  biquadratic tensor.

 We first present an attainable bound for a  biquadratic tensor.

Let $\T = (t_{ijpq}) \in \Re^{d_1 \times d_2 \times d_1 \times d_2}$ be a biquadratic tensor.   We may unfold $\T$ to a $d_1d_2 \times d_1d_2$ matrix $T = (t_{<ij><pq>})$, where $<ij>$ is regarded as one index $<ij>\equiv (i-1)d_1+j = 1, \cdots, d_1d_2$, and $<pq>$ is regard as another index, $<pq> \equiv (p-1)d_1+qd_2 = 1, \cdots, d_1d_2$.  Since $\T$ is biquadratic, matrix $T$ is symmetric.    Note that even if $\T$ is positive semi-definite, $T$ may not be positive semi-definite.  On the other hand, if $T$ is positive semi-definite, $\T$ is always positive semi-definite. If $\mathcal T$ is constructed by a third order tensor as the previous sections, it can be shown that the corresponding matrix $T$ is indeed positive semi-definite.   We do not go to this detail.

We say that $\T$ is rank-one if there are nonzero $\uu \in \Re^{d_1}$ and $\vv \in \Re^{d_2}$ such that $\T = \uu \otimes \vv \otimes \uu \otimes \vv$.

\begin{theorem} \label{t4}
Suppose that $\T = (t_{ijpq}) \in \Re^{d_1 \times d_2 \times d_1 \times d_2}$ is a biquadratic tensor.  Let the symmetric $d_1d_2 \times d_1d_2$ matrix $T$ be constructed as above.   Then the spectral radius of $T$ is an upper bound of the spectral norm of $\T$.  This upper bound is attained if $\T$ is rank-one.   Thus, this upper bound is attainable even if $\T$ is the contraction of a third order tensor $\A$ with $\A$ itself by (\ref{e9}).
\end{theorem}
{\bf Proof} The spectral radius of the symmetric matrix $T$ can be calculated as follows.
\begin{equation} \label{e15}
\rho(T) = \max \left\{ \left| \s^\top T \s \right|  : \s^\top \s = 1, \s \in \Re^{d_1d_2} \right\}.
\end{equation}
We may fold $\s$ to a $d_1 \times d_2$ matrix $S = (s_{ij})$.  Then
$$\s^\top T \s = \langle \T, S \otimes S \rangle \equiv \sum_{i,p=1}^{d_1} \sum_{j,q=1}^{d_2} t_{ijpq}s_{ij}s_{pq}.$$
On the other hand, let $S = \x \otimes \y$ for $\x^\top \x = \y^\top \y = 1, \x \in \Re^{d_1}, \y \in \Re^{d_2}$.   Then $\x^\top \x = \y^\top \y = 1$ implies the vector $\s$, corresponding the matrix $S$, satisfying $\s^\top \s = 1$.   Compare the maximal problems in (\ref{n3}) and (\ref{e15}).   The feasible region of (\ref{n3}) is a subset of the feasible region of (\ref{e15}).  In the feasible region of (\ref{n3}), the two objective functions are equal.   Thus, the optimal objective function value of (\ref{e15}), i.e., the spectral radius of the symmetric matrix $T$, is an upper bound of the optimal objective function value of (\ref{e15}), i.e., the spectral norm of $\T$.   When $\T$ is rank-one, The feasible regions of (\ref{n3}) and (\ref{e15}) are the same, and the objective function values of (\ref{n3}) and (\ref{e15}) are equal.   Then the upper bound is attained in this case.   If $\A$ is rank-one, then $\T = \T^{(3)}$ formed by (\ref{e9}) is also rank-one.   Thus this upper bound is attainable even if $\T$ is formed by (\ref{e9}).
\qed

{\bf Example 1 (Continued)} In this example, we have
$$T^{(3)} = \left(\begin{matrix} 29 & 8 & -13 & 19 \\ 8 & 6 & -14 & 0 \\ -13 & -14 & 35 & 5 \\ 19 & 0 & 5 & 21 \end{matrix}\right).$$
By calculation, the spectral radius of $T^{(3)}$ is $53.1980$.  Its square root is $7.2937$.  This gives an upper bound for the spectral norm of $\A$.
\qed

As in Corollaries \ref{c1} and \ref{c2}, if we take contraction of the first or the second indices of a third order tensor $\A$, we may get different upper bounds for the spectral norm of $\A$.  Hence, there are totally three upper bounds for the spectral norm of a third order tensor.  For Example 1, the two other upper bounds are $8.2529$ and $7.8874$, which are not better than $7.2937$.     Also, this approach involves the calculation of the spectral radius of a $d_1d_2 \times d_1d_2$ (or $d_1d_3 \times d_1d_3$ or $d_2d_3 \times d_2d_3$) symmetric matrix.  When $d_1, d_2$ and $d_3$ are large, this approach involves the calculation of the spectral radius of a high dimensional symmetric matrix.

We now present a different way to obtain this upper bound.   Consider the contraction of $\A$ with itself on the second and third indices.  This result a matrix $B^{(1)} = \left(b^{(1)}_{ij}\right) \in \Re^{d_1 \times d_1}$, with
\begin{equation} \label{4.16}
b^{(1)}_{ij} = \sum_{k=1}^{d_2} \sum_{l=1}^{d_3} a_{ikl}a_{jkl}.
\end{equation}
Then $B^{(1)}$ is a symmetric matrix.

\begin{theorem} \label{t4.2}
Let $\A \in \Re^{d_1 \times d_2 \times d_3}$ and $B$ be constructed by (\ref{4.16}).   The matrix $B^{(1)}$ is positive semi-definite.  The square root of its spectral radius is an upper bound of the spectral norm of $\A$.   This upper bound is equal to the upper stated in Theorem \ref{t4}, when $\T$ in Theorem \ref{t4} is the contraction of $\A$ with $\A$ itself on its first index.   Thus, this upper bound is also attainable.
\end{theorem}
{\bf Proof} We may unfold $\A = (a_{ijk})$ to a $d_1 \times d_2d_3$ matrix $A^{(1)} = (a_{i<jk>})$, where $<jk>$ is regarded as one index $<jk> \equiv (j-1)d_2+k = 1, \cdots, d_2d_3$.   The spectral norm of matrix $A^{(1)}$ can be calculated as
\begin{equation} \label{4.17}
\left\|A^{(1)}\right\| = \max \left\{ \x^\top A^{(1)} \s : \x^\top \x = \s^\top \s = 1, \x \in \Re^{d_1}, \s \in \Re^{d_2d_3} \right\}.
\end{equation}
Compare the maximal problems in (\ref{n1}) and (\ref{4.17}).   The feasible region of (\ref{n1}) is a subset of  (\ref{4.17}).  In the feasible region of (\ref{n1}), the two objective functions are equal.   Hence, the optimal objective function value of  (\ref{4.17}), i.e., the spectral norm of the matrix $A^{(1)}$, is an upper bound of the optimal objective function value of
(\ref{n1}), i.e., the spectral norm of $\A$.   The spectral norm of the matrix $A^{(1)}$ is the largest singular value of $A^{(1)}$, which is equal to the square root of the spectral radius of $A^{(1)}\left(A^{(1)}\right)^\top$.   We now can recognize that $B^{(1)} = A^{(1)}\left(A^{(1)}\right)^\top$.  Thus, $B^{(1)}$ is symmetric and positive semi-definite, and the  square root of its spectral radius is an upper bound of the spectral norm of $\A$.

When $\T = \T^{(1)}$ in Corollary \ref{c1} is the contraction of $\A$ with $\A$ itself on its first index, the upper bound obtained there is equal to the upper bound obtained here.   In fact, in this case, the upper bounds stated in Corollary \ref{c1}, when $\T = \T^{(1)}$, is the square root of the spectral radius of $\left(A^{(1)}\right)^\top A^{(1)}$, while the upper bound given here is the square root of the spectral radius of $A^{(1)}\left(A^{(1)}\right)^\top$.   By linear algebra, they are equal.   Hence, this upper bound is also attainable.
\qed

As $B^{(1)}$ is a $d_1 \times d_1$ symmetric matrix, this approach is relatively easy to be handled.  We may also consider the contraction of $\A$ with itself on the first and third indices, or on the first and second indices.  This results in another way to calculate the two other upper bounds for the spectral norm of $\A$.

\section{Lower Bounds}

We present two attainable lower bounds for the spectral norm of the biquadratic tensor $\T$ in this section.

Let $\T = (t_{ijpq}) \in \Re^{d_1 \times d_2 \times d_1 \times d_2}$ be a biquadratic tensor.  We say that $\T$ is diagonal with respect to its first and third indices if $t_{ijpq} = 0$ whenever $i \not = p$.  We say that $\T$ is diagonal with respect to its second and fourth indices if $t_{ijpq} = 0$ whenever $j \not = q$.

\begin{theorem} \label{t5}
Let $\T = (t_{ijpq}) \in \Re^{d_1 \times d_2 \times d_1 \times d_2}$ be a biquadratic tensor.   A lower bound for the spectral norm of $\T$ is the maximum of the spectral radii of $d_2$ symmetric $d_1 \times d_1$ matrices $(t_{ijpj})$, where $j$ is fixed, for $j = 1, \cdots, d_2$.   This lower bound is attained if $\T$ is diagonal with respect to its second and fourth indices. Another lower bound for the spectral norm of $\T$ is the maximum of the spectral radii of $d_1$ symmetric $d_2 \times d_2$ matrices $(t_{ijiq})$, where $i$ is fixed, for $i = 1, \cdots, d_1$.   This lower bound is attained if $\T$ is diagonal with respect to its first and third indices.
\end{theorem}
{\bf Proof}   Fix $j$.   Let $\y$ be a unit vector in $\Re^{d_2}$ such that its $j$th component is $1$ and its other components are zero.   Then the objective function of (\ref{e9}) is equal to
$$\langle \T, \x \otimes \y \otimes \x \otimes \y \rangle = \sum_{i,p=1}^{d_1}t_{ijpj}x_ix_p.$$
Let $\x$ be the eigenvector of the symmetric matrix $(t_{ijpj})$ such that
$$\sum_{i,p=1}^{d_1}t_{ijpj}x_ix_p = \rho(t_{ijpj}),$$
where $\rho(t_{ijpj})$ is the spectral radius of the symmetric $d_1 \times d_1$ matrices $(t_{\cdot j\cdot j})$.  This is true for $j = 1, \cdots, d_2$.   Hence, the maximum of the spectral radii of $d_2$ symmetric $d_1 \times d_1$ matrices $(t_{\cdot j\cdot j})$, where $j$ is fixed, for $j = 1, \cdots, d_2$, is a lower bound for the spectral norm of $\T$.   Let $\T$ is diagonal with respect to its second and fourth indices.   Then the objective function value of (\ref{e9}) is equal to a convex combination of the spectral radii of $d_2$ symmetric $d_1 \times d_1$ matrices $(t_{\cdot j\cdot j})$, where $j$ is fixed, for $j = 1, \cdots, d_2$.  Then this lower bound is attained in this case.   The other conclusion can be proved similarly.
\qed

{\bf Example 1 (Continued)}  In this example, fix $j = 1$ and $j = 2$, respectively, we have two symmetric matrices
$$\left(\begin{matrix} 29 & - 13 \\ -13 & 35 \end{matrix}\right), \ \ \left(\begin{matrix} 6 & 0 \\ 0 & 21 \end{matrix}\right).$$
Their spectral radii are $45.3417$ and $21$, respectively.  The maximum of these two spectral radii is $45.3417$.  This gives a lower bound of the spectral norm of $\T^{(3)}$.  Its square root is $6.7336$. This gives a lower bound for the spectral norm of $\A$.

Similarly, fix $i = 1$ and $i = 2$, respectively, we have two symmetric matrices
$$\left(\begin{matrix} 29 & 8 \\ 8 & 6 \end{matrix}\right), \ \ \left(\begin{matrix} 35 & 5 \\ 5 & 21 \end{matrix}\right).$$
Their spectral radii are $31.5089$ and $36.6023$, respectively.  The maximum of these two spectral radii is $36.6023$.  Its square root is $6.0500$. This gives another lower bound for the spectral norm of $\A$.
\qed

A question is for which kind of third order tensor $\A$, these two lower bounds are attained.

As in Corollaries \ref{c1} and \ref{c2}, if we take contraction of the first or the second indices of a third order tensor $\A$, we may get different lower bounds for the spectral norm of $\A$.  Hence, there are totally six lower bounds for the spectral norm of a third order tensor.   In particular, for the example in Example 1, if we take contraction of the first index of the third order tensor $\A$ in that example, we get a lower bound $6.7336$ for the spectral norm of $\A$.  As the spectral norm of $\A$ is $6.7673$, this lower bound is about $0.5 \%$ close to the true value.    Surprised by this accuracy, we calculate $1000$ randomly generated examples of $2 \times 2 \times 3$ tensors.  We found that the lower bounds obtained in this way fall within $0.01\%, 0.02\%, 0.05\%, 0.1\%, 0.2\%. 0.5\%, 1\%, 2\%, 5\%, 10\%, 20\%$ and $50\%$ are $4.60\%, 6.40\%, 9.50\%, 13.30\%, 18.80\%, 29.80\%, 40.80\%, 56.00\%, 80.00\%, 94.2\%, 99.50\%$ and $100\%$, respectively.   This shows that for such a third order tensor, there is a big chance to give a good lower bound in this way.

In this approach, spectral radii of $d_i \times d_i$ symmetric matrices for $i = 1, 2, 3$, are calculated.  This only involves relatively low dimensional matrices.  Therefore, this approach is relatively efficient.

\section{Relation}

The first lower bound of $\| \T\|$ in Theorem \ref{t5} may be denoted as
$$L = \max \left\{ \rho((t_{ijpj})) : j \ {\rm is\ fixed},\  j = 1, \cdots, d_2 \right\}.$$
Suppose that $\T = \T^{(3)}$ is constructed by (\ref{e9}) from a third order tensor $\A = (a_{ijk})$.  Then
$$L = \max \left\{ \max \left\{ \sum_{i, p=1}^{d_1} \sum_{k=1}^{d_3} a_{ijk}a_{pjk}x_ix_p : \x^\top \x = 1, \x \in \Re^{d_1} \right\} : j = 1, \cdots, d_2 \right\}.$$
On the other hand, the spectral radius of the matrix $B$, constructed by (\ref{4.16}), is as follows.
$$\rho\left(B^{(1)}\right) = \max \left\{ \sum_{i, p=1}^{d_1} \sum_{j=1}^{d_2} \sum_{k=1}^{d_3} a_{ijk}a_{pjk}x_ix_p : \x^\top \x = 1, \x \in \Re^{d_1} \right\}.$$
Then we find that
$$\rho\left(B^{(1)}\right) \le d_2 L.$$
Combining this with Theorems \ref{t4.2} and \ref{t5}, we have the following theorem.
\begin{theorem} \label{t6}
Let $\A \in \Re^{d_1 \times d_2 \times d_3}$, $L$ and $B^{(1)}$ be constructed as above.    Then we have
$${1 \over d_2}\rho\left(B^{(1)}\right) \le L \le \|\A\| \le \rho\left(B^{(1)}\right) \le d_2L.$$
\end{theorem}

This establishes a range of $\|\A\|$ by either $L$ or $\rho\left(B^{(1)}\right)$.  We may contract on other indices and obtain similar results.   Combining them together, we may get a better range of $\|\A\|$.

\section{Lower Bounds for Nuclear Norms}

Suppose that $\T = (t_{ijpq}) \in \Re^{d_1 \times d_2 \times d_1 \times d_2}$ is a biquadratic tensor. Let the $d_1d_2 \times d_1d_2$ symmetric matrix $T$ be constructed as in Section 4.   Then $T$ is a matrix flattening of the tensor $\T$.   As Lemma 3.1 of \cite{Hu15}, there is a one to one correspondence between the $d_1d_2 \times d_1d_2$ symmetric matrices and $d_1 \times d_2 \times d_1 \times d_2$ biquadratic tensors.  Hence, with an argument similar to the proof of Proposition 4.1 of \cite{Hu15}, we have the following result.

\begin{theorem} \label{t7}
Suppose that $\T = (t_{ijpq}) \in \Re^{d_1 \times d_2 \times d_1 \times d_2}$ is a biquadratic tensor. Let the $d_1d_2 \times d_1d_2$ symmetric matrix $T$ be constructed as in Section 4.  Then $\| T \|_* \le \| \T \|_*$.
\end{theorem}

Combining Theorems 3.3 and \ref{t7}, we have a lower bound for the nuclear norm of a third order tensor, by the nuclear norm of a matrix.    Note that the nuclear norm of a tensor is NP-hard to compute, while the nuclear norm of a matrix is relatively easy to be computed.

Suppose that $\A = (a_{ijk}) \in \Re^{d_1 \times d_2 \times d_3}$, $\T = \T^{(1)}$ is constructed by contraction of $\A$ with $\A$ itself on the first index, and $T^{(1)}$ be the $d_2d_3 \times d_2d_3$ matrix flattening of $\T^{(1)}$.
Then, the square root of the nuclear norm of $T^{(1)}$ also gives a lower bound of the nuclear norm of $\A$.  Let $A^{(1)}$ be the matrix flattening of $\A$ as given in the proof of Theorem \ref{t4.2}.  By \cite{Hu15}, $\left\| A^{(1)}\right\|_*$ also gives a lower bound of $\|\A \|_*$. By the definition of nuclear norms, we find that $\left\|T^{(1)}\right\|_* = \left\| A^{(1)}\right\|_*^2$.    Thus, the lower bound given here is the square root of the lower bound given in \cite{Hu15} for $\|\A\|_*$.  Let $B^{(1)}$ be the $d_1 \times d_1$ symmetric matrix constructed as in Theorem \ref{t4.2}.   With an argument similar to the proof of Theorem \ref{t4.2} and by using the definition of nuclear norms, we may show that
$$\left\| T^{(1)}\right\|_* = \left\| B^{(1)}\right\|_* = \left\| A^{(1)}\right\|_*^2.$$
Since $B^{(1)}$ is symmetric and its dimension is lower, the approach using $B^{(1)}$ may be better than the approach using $A^{(1)}$ in \cite{Hu15}. In \cite{Hu15}, a range of $\|\A\|_*$ is given by $\left\| A^{(1)} \right\|_*$ as:
$$\left\| A^{(1)} \right\|_* \le \| \A \|_* \le \sqrt{ \min \{ d_2, d_3\}} \left\| A^{(1)} \right\|_*.$$
Then we have
$$\sqrt{\left\| B^{(1)} \right\|_*} \le \| \A \|_* \le \sqrt{\min \{ d_2, d_3\} \left\| B^{(1)} \right\|_*}.$$

\section{Final Remarks}

In \cite{JYZ17}, it was shown that the spectral norm and the nuclear norm of a tensor is equal to the spectral norm and the nuclear norm of the Tucker core of that tensor.   As the size of the Tucker core may be smaller than the size of the original tensor, maybe we may combine our results with that approach.

We may also explore more algorithms like that one in \cite{WQZ09} to compute the largest M-eigenvalue of a  positive semi-definite biquadratic tensor, and use them for computing the spectral norm of a third order tensor.

We hope that some further research may explore more applications of the equality between singular values of a third order tensor and M-eigenvalues of the related  positive semi-definite biquadratic tensor.

\bigskip

{\bf Acknowledgments}   The authors are thankful to
Yannan Chen for the discussion on Theorems 4.2 and 7.1, and his calculation, to Yiju Wang and Xinzhen Zhang for their comments, and to Qun Wang for her calculation.

\end{document}